\documentclass[12pt]{amsart}

\usepackage{amsmath,amsthm,amscd,euscript}
\def \r{\mathbb R}
\def \z{\mathbb Z}
\newcommand{\ka}{\kappa}
\newcommand{\e}{\varepsilon}
\newcommand{\la}{\lambda}
\newcommand{\de}{\delta}

\newtheorem{theorem}{Theorem}[section]
\newtheorem{lemma}{Lemma}[section]
\newtheorem{corollary}{Corollary}[section]
\newtheorem{proposition}{Proposition}[section]
\theoremstyle{remark}
\newtheorem{remark}{Remark}[section]
\newtheorem{definition}{Definition}[section]

\numberwithin{equation}{section}

\title
{Energy of a knot: variational principles; Mm-energy.}
\author{Oleg Karpenkov}
\date{November 3th, 2004}
\thanks{Partially supported by INTAS-00-0259,
SS-1972.2003.1 and RBRF-01-01-00660 projects}
\email[Oleg Karpenkov]{karpenk@mccme.ru}
\begin{document}
\input epsf
\maketitle
\tableofcontents
%\vspace{1in}

\section{Introduction}

Let $S^1=\r / (2\pi \z)$ be the circle and $\tau:S^1 \longrightarrow \r ^3$ be
a smooth knot. We will assume that $\tau(t)$ is the arc length
parametrization.
Denote by $D(t_1,t_2)$ the length of the minimal subarc between $t_1$ and
$t_2$ on the circle. Let $|*|$ denote the absolute value of vectors
in $\r ^3$.

Following~\cite{circle}, we denote by
$$
E(\tau)=E_f(\tau)=\iint \limits_{S^1 \times S^1}f(|\tau(t_1)-\tau(t_2)|,
D(t_1,t_2))dt_1dt_2
$$
the energy of the knot $\tau$, where $f(\rho, \alpha)$ satisfies the
following conditions:

1) $f(\rho, \alpha) \in C^{1,1}(U)$, where $U=\{(\rho,\alpha)|
0< \rho \le \alpha, \alpha\le \pi\}$;

2) there exist the following limits:
$$
\lim \limits_
{\genfrac{}{}{0pt}{}{(\rho,\alpha)\in U}{\rho\to 0,\rho/\alpha\to 1}}
f(\rho,\alpha), \quad
\lim \limits_
{\genfrac{}{}{0pt}{}{(\rho,\alpha)\in U}{\rho\to 0,\rho/\alpha\to 1}}
\frac{\partial f(\rho, \alpha)}{\partial \rho}, \quad
\lim \limits_
{\genfrac{}{}{0pt}{}{(\rho,\alpha)\in U}{\rho\to 0,\rho/\alpha\to 1}}
\frac{\partial f(\rho, \alpha)}{\partial \rho}. \quad
$$

Almost all energies are not homothety
invariant, so we will consider only knots of length $2\pi$.

The energy of a knot is not an invariant of the topological class of this knot.
If we make a smooth perturbation of a knot, its energy smoothly changes.
We will consider energies with the following important properties.
The energy is always positive. When a knot crossing tends to a double point,
the energy tends to infinity. So every topological class of knots has a
representative with the minimal value of energy. This knot is called a
{\it normal form} of the class.
It is unknown whether each class has a
unique normal form or not, i.e., whether the normal form for some energy is
an invariant of the topological class or not.
The normal forms satisfy the variational equations considered below.

Some energies have a physical meaning. For example $f=1/(|\tau(t_1)-
\tau(t_2)|)$ is the energy of a charged knot.
Unfortunately, this energy is always
infinite. As long as the charged knot does not break there
must be some other forces
which save the knot. Let us consider a model of such a restriction:
$$
f=\frac{(D^2(t_1,t_2))}{(\tau(t_1)-|\tau(t_2)|)}.
$$
For this energy we will develop our variational principles.

The study of knot energies began with the work of Moffatt (1969)~\cite{Moff1},
and was developed by him in ~\cite{Moff2} following Arnold's
work ~\cite{Arn1}.
The first steps in studying properties of the energies of knots
were made by O'Hara~\cite{O-H1,O-H2,O-H3} and
the first variational principles for polygons in space
were studied by Fukuhara~\cite{pol}.

The aim of this article is to prove that any extremal knot $\tau$
satisfies certain variational equations. The paper is organized as
follows. We start in Section 2 with the definitions and
formulations of the main theorem. In Section 3 we prove this
theorem. In Section 4 we prove that the circle unknot always
satisfies our extremal conditions. Unfortunately the integrals in
the equations do not converge for all possible energies. For
example, they do not converge in the case of the most famous
energy: M\"obius energy.  We discuss this also in Section 4.
Section 5 seems to be independent from the previous sections.  In
Section 5 we represent Mm-energy. The definition of this energy
differs with one regarded above. Nevertheless besides its own
properties Mm-energy has some similar with M\"obius energy
properties.

This work is partially published (see~\cite{Kar1} and~\cite{Kar2}).

The author is grateful to professor A.~B.~Sossinsky for
constant attention to this work.

\vspace{5mm}
\section{Notation and definitions}

Mostly we will work with knots of fixed length $2\pi$.
So let $S^1=\r / (2\pi \z)$ be the circle and let
$\tau:S^1 \longrightarrow \r ^3$ denote some smooth knot of length $2\pi$.
Let $\tau(t)$ be the arc length
parametrization.

By $\ka (t)$ we denote the curvature
at $t$ and $R(t)=1/ \ka (t)$, the radius of curvature at $t$.

\begin{definition}\label{d-1}
Given a smooth knot $\tau:S^1 \longrightarrow \r ^3$ and a point
$t_0 \in S^1$, a {\it locally perturbed knot} is a knot (denoted
by $\tau_{t_0,\e}$) such that

a) $|\tau(t)-\tau_{t_0,\e}(t)|<\e ^2$ if $D(t_0,t) \le \e$
and $\tau(t)=\tau_{t_0,\e}(t)$ if $D(t_0,t)> \e$;

b) $|\ka (t)-\ka_{t_0,\e}(t)|<\e$ for $D(t_0,t)< \e$;

c) $\tau_{t_0,\e}(t_0+\lambda)=\tau_{t_0,\e}(t_0)+\lambda \dot {\tau}_{t_0,\e}
(t_0)+({\lambda ^2}/{2})\ddot {\tau}_{t_0,\e}(t_0)+
o(\e^2)$ if $D(t_0,t_0+\lambda)\le \e$.
\end{definition}

Note that at the points $t_0-\e$ and $t_0+\e$ the curvature is
not restricted.

The length of the knot $\tau_{t_0,\e}$ can change, but
we regard knots of length $2\pi$ only.
One of the ways to solve this problem is to consider the restriction of
the set of locally perturbed knots to the set of knots of constant
length  $2\pi$, but this definition is unsatisfactory.
Indeed, let a knot $\tau$ in some neighborhood of the point $t_0$
be a piece of a straight line. Then the set of locally
perturbed knots at the point $t_0$ of length $2 \pi$ consists of the
knot $\tau$ only.

We will extend this set in the following way.

\begin{definition}\label{d-2}
Let the length of $\tau_{t_0,\e}$ be $(1+\de)2\pi$.
The {\it locally perturbed length $2\pi$ knot}
$\tilde \tau_{t_0,\e}$ is the knot obtained from $\tau_{t_0,\e}$
by homothety with coefficient $1/(1+\de)$ and center at the origin.
We also say that the knot $\tilde \tau$ is {\it associated}
with the knot $\tau$.
\end{definition}

Consider any $\tau_{t_0,\e}$. We will show later that
$\de = c_1\e ^3 +\circ (e^3)$. Thus by Definition~\ref{d-1} we have
$$
|\tau_{t_0,\e}(t_1)-\tau_{t_0,\e}(t_2)|=|\tau(t_1)-\tau(t_2)|+
c_2 (t_1,t_2)\e ^2 +o(e^2)
$$
if $D(t_0,t_1)<\e$ or $D(t_0,t_2)<\e$.
Then we may conclude that
$$
E(\tau_{t_0,\e})=E(\tau)+c_3 \e^3+o(\e^3) \quad \mbox{and} \quad
E(\tilde \tau_{t_0,\e})=E(\tau)+c_4 \e^3+o(\e^3).
$$
The coefficients $c_3$ and $c_4$ of the term $\e^3$ will be called
the {\it variation} and denoted by $Var(\tau_{t_0,\e})$ and
$Var(\tilde \tau_{t_0,\e})$ respectively.

Now all is prepared for the definition of a locally extremal point of a knot.

\begin{definition}\label{d-3}
Any $t_0 \in S^1$ is called {\it locally extremal point} of $\tau$
if $Var(\tilde \tau_{t_0,\e})=0$
for each locally perturbed knot $\tilde \tau_{t_0,\e}$ of length $2\pi$.
\end{definition}

\begin{definition}\label{d-4}
The knot $\tau$ is said to be {\it locally extremal} if all its points
are locally extremal.
\end{definition}

Let us find necessary and sufficient conditions for the point
$t_0$ be locally extremal.
We denote the vector product of two vectors $a$ and $b$
by $[a,b]$. By $(a,b,c)$ we denote the mixed product (oriented volume)
of the vectors $a$, $b$ and $c$.
Let $\dot {\tau}(t)$ be the velocity vector and $\ddot {\tau}(t)$ be the
acceleration vector. Now we define the functions
$\Psi(t_0,t)$ and $\Phi(t_0,t)$.
$$
\begin{array}{c}
\Psi(t_0,t)=
\left\{
\begin{array}{cl}
\Bigl( \frac{\dot{\tau}(t_0)}{|\dot{\tau}(t_0)|},
\frac{\ddot{\tau}(t_0)}{|\ddot{\tau}(t_0)|},
\frac{\tau(t)-\tau(t_0)}{|\tau(t)-\tau(t_0)|} \Bigr)
&, $ if $ \ddot{\tau}(t_0) \ne 0;\\
\Bigl( \frac{\tau(t)-\tau(t_0)}{|\tau(t)-\tau(t_0)|},
\frac{\dot{\tau}(t_0)}{|\dot{\tau}(t_0)|}\Bigr)
&, $ if $ \ddot{\tau}(t_0)=0.
\end{array}
\right.\\
\Phi(t_0,t)=
\left\{
\begin{array}{cl}
\Bigl( \frac{\dot{\tau}(t_0)}{|\dot{\tau}(t_0)|},
\frac{\tau(t)-\tau(t_0)}{|\tau(t)-\tau(t_0)|},
\Bigl[\frac{\dot{\tau}(t_0)}{|\dot{\tau}(t_0)|},
\frac{\ddot{\tau}(t_0)}{|\ddot{\tau}(t_0)|}\Bigr] \Bigr)
&, $ if $ \ddot{\tau}(t_0) \ne 0;\\
0 &, $ if $ \ddot{\tau}(t_0)=0.
\end{array}
\right.
\end{array}
$$
Note that $|\dot{\tau}(t_0)|=1$ and $|\tau(t)-\tau(t_0)| \ne 0$ if $t \ne t_0$.
Thus $\Psi$ and $\Phi$ are well defined.

We also remark that $\Psi(t_0,t)=\sin\psi(t_0,t)$,
where $\psi(t_0,t)$ is the angle between the vector $\tau(t)-\tau(t_0)$
and the oriented plane spanning of $\dot{\tau}(t_0)$ and $\ddot{\tau}(t_0)$.
The function $\Phi$ has a similar representation:
$\Phi(t_0,t)=\sin\phi(t_0,t)$,
where $\phi(t_0,t)$ is the angle between the vector $\tau(t)-\tau(t_0)$
and the oriented plane spanning of $\dot{\tau}(t_0)$ and
$[\dot{\tau}(t_0),\ddot{\tau}(t_0)]$.
(See Fig. 1).
These angles can be either positive or negative.
\begin{figure}
$$\epsfbox{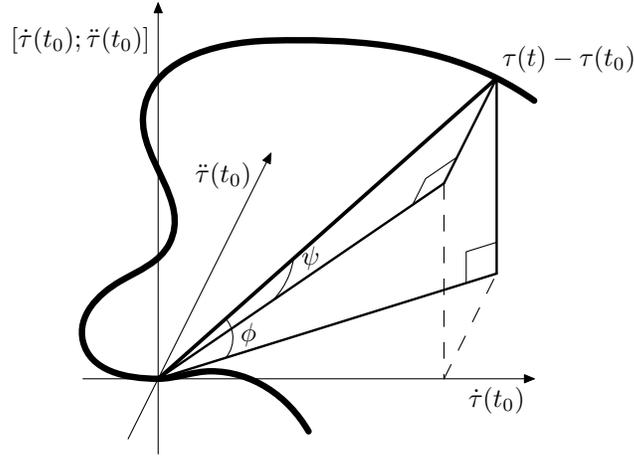}$$
\caption{The geometric interpretation of $\psi(t_0,t)$ and $\phi(t_0,t).$}
\end{figure}

\begin{theorem}\label{t-1}
Let $\tau$ be a smooth knot. The point $t_0$ is a locally extremal point
of $\tau$ if and only if the following conditions hold:
$$
\begin{array}{l}
\begin{aligned}
\displaystyle
V_1(t_0)&:=&\frac{2}{3R(t_0)}\Biggl(
4 \int \limits_{S^1} \Bigl( f+ R(t_0)\Phi(t_0,t)
\frac{\partial f}{\partial \rho} \Bigr) dt-
\frac{1}{\pi}\iint \limits_{S^1 \times S^1} \Bigl( 2f+D(t_1,t_2)
\frac{\partial f}{\partial \rho}+\\
\displaystyle
&&|\tau (t_1)-\tau (t_2)| \frac{\partial f}{\partial \alpha} \Bigr) dt_1dt_2+
2\iint \limits_{A} \frac{\partial f}{\partial \alpha} dt_1dt_2 \Biggr)=0;
\end{aligned}\\
\displaystyle
V_2(t_0):=\frac{4}{3R(t_0)} \int \limits_{S} \frac{\partial f}{\partial \rho}
\Psi(t_0,t) dt=0.
\end{array}
$$
Here $A \subset S^1 \times S^1$ is the set of points $(t_1, t_2)$ such that
$D(t_1,t_2)=D(t_1,t_0)+D(t_0,t_2). $
\end{theorem}

\begin{corollary}\label{c-1}
A knot $\tau$ is locally extremal if and only if almost all of its
points are locally extremal, i.e.,
$$
\int \limits_{S^1} \Bigl( V_1^2(t)+V_2^2(t) \Bigr) dt=0.
$$
\end{corollary}

\vspace{5mm}
\section{Proofs}

Let $t_0$ be any point of $S^1$. We choose orthonormal
coordinates in $\r ^3$ such that $\tau(t_0)$ is on the $(X,Y)$-plane,
$\tau(t_0-\e)$ and $\tau(t_0+\e)$ lie symmetrically on the $X$-axis.
If $\tau(t_0-\e)$, $\tau(t_0)$ and $\tau(t_0 + \e)$ are on the same
line, then we make any possible choice of the $Y$-axis.
Finally, we choose the $Z$-axis such that the orientation of the
$(X,Y,Z)$-space is positive (see Fig. 2a)).

Let $P_{\e}$ be the class of parabolic arcs and one segment such that all
the parabolas have their vertex in the $(Y,Z)$-plane,
$\tau(t_0-\e)$ and $\tau(t_0+\e)$ are the endpoints of the arcs,
and the endpoints of the segment are $\tau(t_0-\e)$ and $\tau(t_0 + \e)$.
Each parabola can be specified by two parameters $(\lambda, \gamma)$,
where $2\lambda$ is the ``acceleration" and $\gamma$ is the angle between
the $(X,Y)$-plane and the plane containing the parabola (see Fig. 2b)).
Notice also that $(0,\gamma)$ is some segment.

\begin{figure}
$$\epsfbox{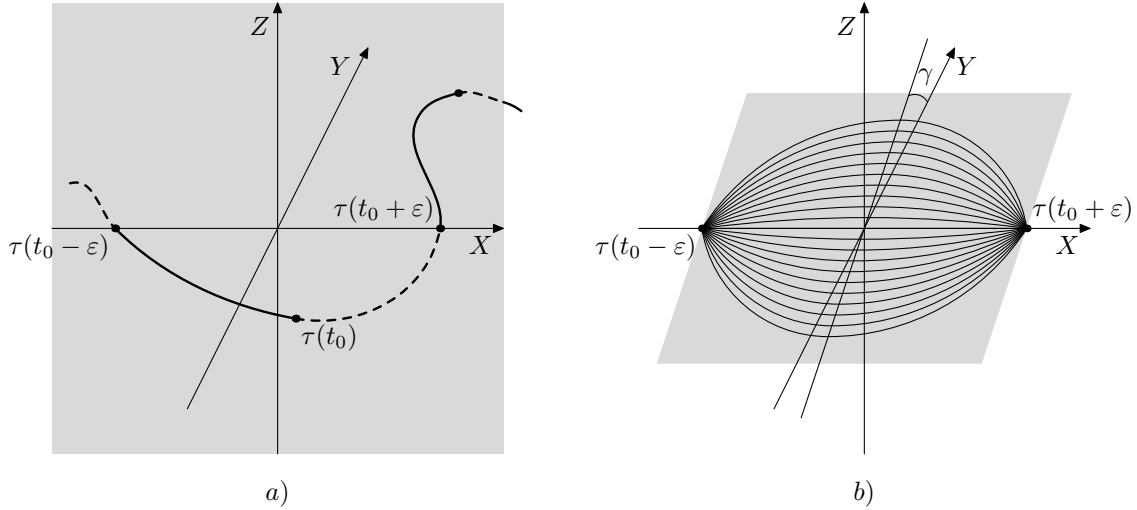}$$
\caption{a)The choice of $X$, $Y$ and $Z$-axes. b)The parabolic arcs
$(\alpha,\gamma)$ where $\gamma$ is fixed.}
\end{figure}

Denote by $M_{P,t_0,\e}$ the $2$-dimensional
set of knots $\tau_{t_0,\e, \la, \gamma}$, where the curve
connecting $\tau(t_0-\e)$ and $\tau(t_0+\e)$ belongs to the class $P_{\e}$
with the following property: the knot
$(\tau_{t_0, \e, \la, \gamma}+ \tau)/2$ is a locally perturbed
knot.
Denote by $\tilde M_{P,t_0,\e}$ the set of knots associated with
the knots in the class $P_{\e}$.

\begin{theorem}\label{t-2}
Let $\tau$ be a smooth knot. The point $t_0$ is a locally extremal point
if and only if $Var(\tilde {\tau}_{t_0,\e})=0$
for each locally perturbed $($at $t_0$$)$ knot
$\tilde \tau_{t_0,\e} \in \tilde M_{P,t_0,\e} $.
\end{theorem}

{\bf Proof of Theorem~\ref{t-2}}.

\begin{figure}
$$\epsfbox{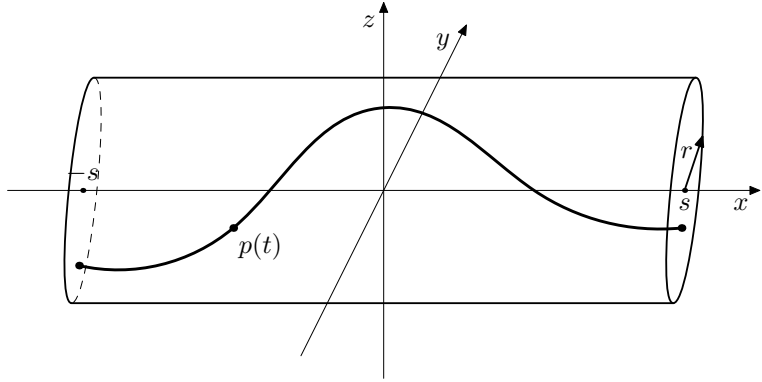}$$
\caption{The trajectory of the point $p(t)$ inside the cylinder $C$.}
\end{figure}

We begin the proof with the following lemma.

\begin{lemma}\label{l-2.1}
Let $C=\{(x,y,z) \in \r ^3 |\sqrt{y^2+z^2}<r,|x|<s\}$ be
a cylinder. Suppose a point moves inside $C$ with velocity of constant
modulus $1$ and so that the absolute value of its acceleration
is bounded by $K$ $($see Fig.3$)$. Let $x(0)=-s$, $x(T)=s$, $s \gg r$
and $K<1/(4r)$. Then the length of the trajectory of a point
$($i.e. T\/$)$ is bounded:
$$
T<\frac{2s}{\sqrt{1-4Kr}}.
$$
\end{lemma}

First let us prove that $\dot y^2 (t_0) < 2Kr$.
We first consider the case for which $x(t_0)<0$,
$\dot x (t_0)>0$ and $y(t_0)>0$; then
$$
y(t)=y(t_0)+\int \limits_{t_0}^{t} \dot y (\xi) d\xi<r.
$$
By the assumption, we have
$$
\dot y(\xi) =\dot y (t_0)+\int \limits_{t_0}^{\xi} \ddot y (\zeta)d\zeta>
\dot y (t_0)-\int \limits_{t_0}^{\xi} Kd\zeta=
\dot y (t_0)-(\xi-t_0)K.
$$
It follows that
$$
y(t)>y(t_0)+\int \limits_{t_0}^{t} \dot y (x_0)-(\xi-t_0)Kd\xi=
y(t_0)+(t-t_0) \dot y (t_0) - \frac{(t-t_0)^2}{2}K.
$$
But $y(t_0)>-r$ and $y(t)<r$, so
$$
(t-t_0) \dot y (t_0) - \frac{(t-t_0)^2}{2}K-2r<0.
$$
By assumption $x<0$ and $s \gg r$, so the vertex of the parabola is
at the point $t-t_0={\dot  y (t_0)}/{K}<s$.
This yields the inequality $\dot y ^2 (t_0) < 2Kr$.

The proof for the cases in which
$\dot x(t_0)>0$ and $y(t_0)<0$;
$\dot x(t_0)<0$ and $y(t_0)>0$;
$\dot x(t_0)<0$ and $y(t_0)<0$ is similar.

Secondly, we claim that  $\dot z ^2(t_0) < 2Kr$. The proof is similar
to the inequality for $\dot y ^2(t_0)$.

By the previous statements, it follows that
$$
\dot x ^2(t_0) =1-\dot y ^2(t_0)-\dot z ^2(t_0)>1-4Kr>0
$$
for every
$t_0 \in [0,T]$. So we have $T<2s/(1-4Kr)$

This completes the proof of Lemma~\ref{l-2.1}

We continue the proof with a generalization of the previous lemma.

\begin{figure}
$$\epsfbox{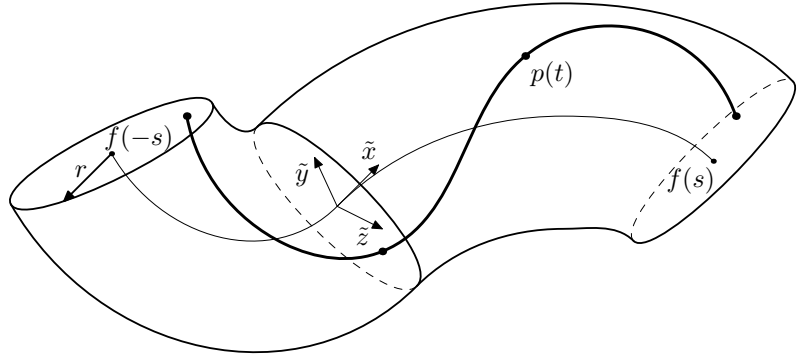}$$
\caption{The trajectory of the point p(t) inside the cylinder $C_f$.}
\end{figure}

\begin{lemma}\label{l-2.2}
Let $f:[-s,s] \longmapsto \r ^3$ be a unit-length smooth map,
let the curvature of $f$ be bounded $(|\ddot f(t)|<K_1)$ and
$sK_1 < 1$.
Let $D^2(t) \in \r ^3$, where $t \in [-s,s]$ is the disk of radius
$r$ centered at $f(t)$
with the plane of the disc orthogonal to $\dot f$. Let also $rK_1<1$.
Denote by $C_f=\bigcup_{[-s,s]} D^2(t)$ the tubular neighborhood of
the curve $f$.
Suppose a point $p(t)=(x(t),y(t),z(t))$ moves inside $C_f$ $($see Fig. 4$)$ with
velocity of constant absolute value $1$ and let the absolute value of
its acceleration be bounded by $K_2$.
Let $p(0) \in D^2(-s)$, $p(T) \in D^2(s)$. Let
$s \gg r$ and
$$
K_2+\frac{1}{1-rK_1}K_1<\frac{1}{4r}.
$$
Then the length of the trajectory of the point $($i.e., T$)$ is bounded and
$$
2s(1-rK_1)<T<\frac{2s(1+rK_1)}{\sqrt{1-4K_2r}}.
$$
\end{lemma}

Let us define $\tilde x =t$.

Now we describe some map $\pi$ from $C_f$ to the standard
cylinder $C$ (see Fig. 3). Let
$$
\pi (D^2(\tilde x))=\{(\tilde x,y,z) \in \r | \sqrt{y^2+z^2}\}
$$
be isometric images of the disk $D^2(\tilde x)$ for each $\tilde x\in[-s,s]$.
If we fix a preimage $\tilde y$-axis of the $y$-axis and
a preimage $\tilde z$-axis of the $z$-axis in the disc $D^2(\tilde x)$
for each $\tilde x \in [-s,s]$, then the map will be completely described.
As long as $sK_1<1$ and $rK_1<1$, this map is well defined and
the manifold $N_f=\bigcup_{[-s,s]} \partial D^2(t)$ with boundary
$\partial D^2(-s)\cup \partial D^2(s)$ is smooth.

Let $\pi (\tilde y_{-s})=(-s,r,0)$ for some $\tilde y_{-s} \in \partial
D^2(-s)$. Consider the vector field on $N_f$ with the following
property: if the point $q$ lies on the circle $\partial D^2(\tilde x)$,
then the vector $v_q$ equals $\dot f(\tilde x)$; this means that
$v_q$ is the unit-length vector orthogonal to the disc $D^2(\tilde x)$ with the
corresponding direction. Denote the integral trajectory of this field
passing through the point $\tilde y_{-s}$ by
$\tilde y=\{\tilde y(\tilde x)| \tilde x \in [-s,s]\}$.
This trajectory defines the $\tilde y$ coordinate
in each disc $D^2(\tilde x)$. Finally we define the unit-length
$\tilde z$-vector as
the vector product of the unit-length $\tilde x$-vector and unit-length
$\tilde y$-vector (in each $D^2(\tilde x)$).

The image $\pi(p)$ of the point $p$ moves inside $C$. We denote $\pi (p)$
by $\hat p$.
Notice that
$$
\frac{|\dot p(t)|}{|\dot{\hat p}(t)|}=
\frac{1}{|\dot{\hat p}(t)|} \in [1-rK_1,1+rK_1].
$$
Note also that if the curvature of the trajectory is $K$
at some point $p(t)$,
then the curvature of the image of this trajectory will be
$$
\hat K < K+\frac{1}{\frac{1}{K_1}-r}
$$
at the point $\hat p(t)$, as can be easily shown.

Now Lemma~\ref{l-2.2} follows from Lemma~\ref{l-2.1}.

\vspace{0.2cm}
We continue the proof of Theorem~\ref{t-2}.
Let $\tau_{t_0,\e}$ be any locally perturbed knot at the point $t_0$
and let $t \in S^1$ such that $D(t_0,t)<\e$. Consider
$$
P_{\tau_{t_0,\e}}(t)= \tau_{t_0,\e}(t_0) +
(t-t_0) \bigl(\dot {\tau}_{t_0, \e}(t_0)+c_1\bigr)+
\frac{(t-t_0)^2}{2}\bigl(\ddot{\tau}(t_0)+c_2\bigr)
$$
We choose the constants $c_1=o(\e ^2)$ and $c_2=o(\e ^2)$ so that
$$
P_{\tau_{t_0,\e}}(t+\e)=\tau_{t_0,\e}(t+\e), \quad
P_{\tau_{t_0,\e}}(t-\e)=\tau_{t_0,\e}(t-\e).
$$
Here we take the unit-length parametrization and denote the length
of curves by $l(*)$.
Then $P_{\tau_{t_0,\e}}(t)$ is a parabolic arc in the $\e$-neighborhood
of the point $t_0$.
From Lemma~\ref{l-2.2} it follows that $|P_{\tau_{t_0,\e}}(t)
-\tau_{t_0,\e}(t)|=o(\e^2)$
and also $l(P_{\tau_{t_0,\e}}(t))=l(\tau_{t_0, \e})+o(\e^3)$.
So $E_{\tau_{t_0,\e}}-E_{P_{\tau_{t_0,\e}}(t)}=o(\e^3).$

Consider the perturbed curve $\tau_{t_0,\e,\la,\gamma}$ passing through
the point $\tau_{t_0,\e}(t_0)$. We have
$$
|\tau_{t_0,\e,\la,\gamma}(t)-\tau_{t_0,\e}(t)|<\e.
$$
We also have $E_{\tau_{t_0,\e,\la,\gamma}}-P_{E_{\tau_{t_0,\e}}}=o(\e^3)$.

Finally we conclude that
$E_{\tau_{t_0,\e,\la,\gamma}}-E_{\tau_{t_0,\e}}=o(\e^3)$.

One can see that the knot $\tau_{t_0,\e,\la,\gamma}$ belongs
$M_{P,t_0,\e}$. We note again that
$l(P_{\tau_{t_0,\e}}(t))=l(\tau_{t_0, \e})+o(\e^3)$.
Hence
$$
E_{\tilde\tau_{t_0,\e,\la,\gamma}}-E_{\tilde \tau_{t_0,\e}}=o(\e^3).
$$
By definition, the knot $\tilde \tau_{t_0,\e,\la,\gamma}$ belongs
$\tilde M_{P,t_0,\e}$.
This completes the proof of Theorem~\ref{t-2}.

{\bf Proof of Theorem~\ref{t-1}}
Without loss of generality, we put
$$
t_0=0, \qquad \gamma=o(1),\quad \mbox{and} \quad \lambda=1/(2R(0))+o(1),
$$
where R(0) is the radius of curvature at the point $0$.
According to Theorem~\ref{t-2}, we can consider only the class $\tilde M_P$
of knots. Let $\tilde \tau_{0, \e, \la, \gamma}$ be a knot in $\tilde K_P$.
Denote
$$
\Delta:=\Bigl[\frac{\e}{1+\delta},\frac{\e}{1+\delta}\Bigr]\subset S^1.
$$
Now note that for any $\tau$ we have

$$
\begin{array}{c}
\displaystyle
E(\tau)=
\iint \limits_{S^1 \times S^1} f dxdy=
2 \iint \limits_{\Delta \times S^1} f dxdy-
\iint \limits_{\Delta \times \Delta} f dxdy+
\iint \limits_{A \setminus(\Delta \times S^1 \bigcup S^1 \times \Delta)}
f dxdy+
\\
\displaystyle
\iint \limits_{S^1 \times S^1 \setminus A} f dxdy =:
2E_1(\tau)-E_2(\tau)+E_3(\tau)+E_4(\tau).
\end{array}
$$

Here $f=f(\rho(\tau(x),\tau(y)),\alpha(\tau(x),\tau(y))).$
Further note that
$$
Var(\tau)=2Var_1(\tau)-Var_2(\tau)+Var_3(\tau)+Var_4(\tau),
$$
where $Var_i$ is the variation of $E_i$.

First we calculate $Var_1$. We recall that $\sin \phi = \Phi$ and
$\sin \psi = \Psi$.
\begin{lemma}\label{l-2}
$$
Var_1(\tilde \tau_{0,\e,\la,\gamma})=
\frac{4}{3}\Bigl( \int \limits_{S^1} \frac{f}{R(0)}
+\sin\phi \frac{\partial f}{\partial\rho}dy\Bigr) \bigl(\lambda-
\frac{1}{2R(0)}\bigr)+ \frac{2}{3R(0)}\Bigl(
\int \limits_{S^1}\sin\psi \frac{\partial f}{\partial\rho} dy\Bigr) \gamma.
$$
\end{lemma}

The length of the arc of the parabola is $2\e+\frac{2}{3}\lambda^2\e^3+
o(e^3)$. So $\delta=\frac{2}{3}\lambda^2\e^3$.
Note also that
the coefficient of homothety is $o(e^2)$ and thus
$Var_1(\tau_{0,\e,\la,\gamma})=Var_1(\tilde \tau_{0,\e,\la,\gamma})$.
Let
$$
(a,b,c)=(a(t),b(t),c(t))=\tau_{0,\e,\la, \gamma}(t), \quad
\ell=\ell (t)=\sqrt{a(t)^2+b(t)^2+c(t)^2}, \quad f=f(\rho,\alpha).
$$
Thus we have
$$
\begin{array}{c}
\displaystyle
E_1(\tau_{0,\e,\la})=\int \limits_{S^1} \int \limits_{-\e}^{\e}
\biggl( \Bigl[ 1+(2\lambda {t_1})^2\Bigr] ^{1/2}\\
\displaystyle
f\Bigl( \bigl[ (({t_1}-a({t_2}))^2+((\lambda {t_1}^2-\lambda \e^2)\cos\gamma
-b({t_1}))^2+
((\lambda {t_1}^2-\lambda \e^2)\sin\gamma -c({t_1}))^2) \bigr] ^{1/2},\\
\displaystyle
D({t_1},{t_2}) \Bigr) \biggr) d{t_1}d{t_2}+ o(\e^3)=\\
\displaystyle
\int \limits_{S^1} \int \limits_{-\e}^{+\e}
\biggl( \Bigl(1+2\lambda^2{t_1}^2+o({t_1}^2)\Bigr)
\Bigl(f+\Bigl(\frac{{t_1}^2-a{t_1}+\lambda (\e^2-{t_1}^2)
(b\cos \gamma+c\sin \gamma)}{\ell}-\frac{a^2{t_1}^2}{2\ell^2}\Bigr)\times
\\
\displaystyle
{\frac{\partial f}{\partial \rho}}+
\frac{a^2{t_1}^2{\frac{\partial ^2  f}{\partial \rho ^2}}}{2\ell}+D(0,{t_1})
{\frac{\partial f}{\partial \alpha}} \Bigr) \biggr) d{t_1}d{t_2}+o(\e^3)=\\
\end{array}
$$
$$
\begin{array}{c}
\displaystyle
=\int \limits_{S^1} \int \limits_{-\e}^{+\e}
\biggl( 2\lambda^2{t_1}^2f+f+\Bigl( \frac{{t_1}^2-a{t_1}+\lambda
(\e^2-{t_1}^2)(b\cos \gamma+ c\sin \gamma)}{\ell}-\\
\displaystyle
\frac{a^2{t_1}^2}{2\ell^2} \Bigr) {\frac{\partial f}{\partial \rho}}+
\frac{a^2{t_1}^2{\frac{\partial ^2  f}{\partial \rho ^2}}}{2\ell}+D({t_1},0)
{\frac{\partial f}{\partial \alpha}} \biggr) d{t_1}d{t_2}+o(\e^3)=\\
\displaystyle
=\int \limits_{S^1}
\biggl(\frac{2\lambda^2\e^3}{3}f+\e f+\Bigl( \frac{2\e^3+
4\lambda \e^3(b\cos \gamma+ c\sin \gamma)}{3\ell}-\frac{a^2\e^3}{2\ell^2}\Bigr)
{\frac{\partial f}{\partial \rho}}+\frac{a^2\e^3
{\frac{\partial ^2  f}{\partial \rho ^2}}}{2\ell} \biggr) d{t_2}+\\
\displaystyle
\int \limits_{S^1} \int \limits_{-\e}^{+\e} D({t_1},0)
{\frac{\partial f}{\partial \alpha}} d{t_1}d{t_2}+o(\e^3)
\end{array}
$$
This yields
$$
\begin{array}{c}
\displaystyle
dE_1(\la,\gamma)=
d\biggl( \int \limits_{S^1}
\Bigl(\frac{2\lambda^2\e^3}{3}f+\Bigl(\frac{4\lambda \e^3(b\cos \gamma+
c\sin \gamma)}{3\ell}\Bigr) {\frac{\partial f}{\partial \rho}} \Bigr) d{t_2}
+o(\e^3)\biggr)=\\
\displaystyle
\biggl( \int \limits_{S^1}
\Bigl(\frac{4\lambda \e^3}{3}f+\Bigl(\frac{4 \e^3(b\cos \gamma+
c\sin \gamma)}{3\ell}\Bigr){\frac{\partial f}{\partial \rho}}\Bigr) d{t_2}
+o(\e^3)\biggr)d\lambda+\\
\displaystyle
\biggl( \int \limits_{S^1}
\Bigr(\frac{4 \e^3\lambda (-b\sin \gamma+ c\cos \gamma)}{3\ell}\Bigl)
{\frac{\partial f}{\partial \rho}} d{t_2}+o(\e^3)\biggr)d\gamma.
\end{array}
$$
Finally  we substitute
$$
\frac{b}{\ell}=\sin\phi, \quad \frac{c}{\ell}=\sin\psi, \quad
\gamma=o(1), \quad \lambda=\frac{1}{2R(0)}+o(1),
$$
where $R(0)$ is the radius of curvature
at the point $0$, obtaining
$$
Var_1(\tilde \tau_{0,\e,\la,\gamma})=
\frac{4}{3}\Bigl( \int \limits_{S^1}
\frac{f}{R(0)}+\sin\phi {\frac{\partial f}{\partial \rho}}d{t_2}\Bigr)
\Bigl(\lambda-\frac{1}{2R(0)}\Bigr)+
\frac{2}{3R(0)}\Bigl( \int \limits_{S^1}\sin\psi {\frac{\partial f}
{\partial \rho}}d{t_2}\Bigr) \gamma.
$$
The proof of Lemma~\ref{l-2} is complete.

\begin{lemma}\label{l-3} $Var_2=0.$
\end{lemma}

Since $E_2(\tau)-E_2(\tilde \tau_{0,\e})=o(\e^3)$,
we immediately have $Var_2=0$.

\begin{lemma}\label{l-4}
$$
Var_3=\Bigl( \frac{2}{3\pi R(0)}\iint
\limits_{A \setminus(\Delta \times S^1 \bigcup S^1 \times \Delta)}
-2f-\ell {\frac{\partial f}{\partial \rho}}+(2\pi-D({t_1},{t_2}))f_{\la}
d{t_1}d{t_2} \Bigr) \bigl(\la- \frac{1}{2R(0)}\bigr).
$$
\end{lemma}

The following calculations prove this lemma.

$$
\begin{array}{c}
\displaystyle
E_3(\tau_{0,\e,\la})=\iint \limits_{A \setminus(\Delta \times S^1 \bigcup S^1 \times \Delta)}
fd\tilde {t_1} d \tilde {t_2}=\\
\displaystyle
\iint \limits_{A \setminus(\Delta \times S^1 \bigcup S^1 \times \Delta)}
f\biggl( \ell \Bigl(1-\frac{2\la^2\e^3}{3\pi}\Bigr),
D\Bigl( \Bigl(|{t_2}-{t_1}|+\frac{4\la^2\e^3}{3}\Bigr)
\Bigr(1-\frac{2\la^2\e^3}{3\pi}\Bigr),0\Bigr) \biggr)\\
\displaystyle
d\Bigl( \bigl(1-\frac{2\la^2\e^3}{3\pi}\bigr){t_1}\Bigr)
d\Bigl( \bigl(1-\frac{2\la^2\e^3}{3\pi}\bigr){t_2}\Bigr)+ o (\e ^3)=\\
\displaystyle
\iint \limits_{A \setminus(\Delta \times S^1 \bigcup S^1 \times \Delta)}
\biggl(f+\Bigl(-\frac{2}{\pi}f-\frac{\ell}{\pi}
\frac{\partial f}{\partial \rho}+ 2f_\la-\frac{D({t_1},{t_2})f_{\la}}{\pi}
\Bigr)\frac{2\la^2 \e^3}{3}\biggr)d{t_1}d{t_2}+ o (\e ^3).
\end{array}
$$
Let us remark that
$$
\iint \limits_{\Delta \times S^1 \bigcup S^1 \times \Delta}
\left(-\frac{2f}{\pi}-\frac{\ell {\frac{\partial f}{\partial \rho}}}{\pi}+
2f_\la-\frac{D({t_1},{t_2})f_{\la}}{\pi}
\right)\frac{2\la^2 \e^3}{3}d{t_1}d{t_2}=o(\e^3).
$$
Therefore
$$
Var_3=\biggl( \frac{2}{3\pi R(0)}\iint
\limits_{A \setminus(\Delta \times S^1 \bigcup S^1 \times \Delta)}
-2f-\ell {\frac{\partial f}{\partial \rho}}+\bigl(2\pi-D({t_1},{t_2})
\bigr)f_{\la}d{t_1}d{t_2} \biggr)\Bigl(\la-\frac{1}{2R(0)} \Bigl).
$$

\begin{lemma}\label{l-5}
$$
Var_4=\biggl( \frac{2}{3\pi R(0)}\iint
\limits_{S^1 \times S^1 \setminus A}
-2f-\ell {\frac{\partial f}{\partial \rho}}-D({t_1},{t_2}){\frac{\partial f}
{\partial \alpha}}d{t_1}d{t_2} \biggl)\Bigr(\la-\frac{1}{2R(0)} \Bigl)
$$
\end{lemma}

The proof of this lemma is similar to the previous one.

Lemmas~\ref{l-2}-\ref{l-5} complete the proof of Theorem~\ref{t-1}.
\vspace{5mm}

\section{Corollaries}

In~\cite{circle} it is shown that the circle is not always the global maximum,
or the global minimum for the energy considered.
Let us show that circle is a locally extremal knot for
any energy $E$ satisfying the conditions 1), 2) of the Introduction.

\begin{corollary}\label{c-2}
The circle is always a locally extremal knot.
\end{corollary}

If $\tau$ is a circle, then
$$
\ell(t_1,t_2)=2\sin\frac{t_2-t_1}{2}, \quad
R(t_1)=1, \quad \psi(t_1,t_2)=0, \quad \phi= \frac{t_2-t_1}{2}.
$$
So $V_2(t_1)=0$ for any $t_1 \in S^1$. Further

$$
\begin{array}{c}
\displaystyle
V_1(t_1)=\frac{1}{3}\biggl(
8\int \limits_{S^1}f+\sin\Bigl(\frac{|t_1|}{2}\Bigr)f_\rho dt_1 -
\frac{2}{\pi}\iint\limits_{S^1 \times S^1} 2f+
2\sin\Bigl(\frac{|t_2-t_1|}{2}\Bigr)\frac{\partial f}{\partial \rho}\\
\displaystyle
+D(t_1,t_2)\frac{\partial f}{\partial \alpha}dt_1dt_2+
4\iint\limits_{A} \frac{\partial f}{\partial \alpha} dt_1dt_2 \biggr) =\\
\displaystyle
\frac{1}{3}\biggl( 8\int \limits_{S^1}f+\sin\Bigl(\frac{|t_1|}{2}\Bigr)
\frac{\partial f}{\partial \rho} dt_1 -4\int \limits_{S^1} 2f+2\sin\Bigl(\frac{|t_1|}{2}\Bigr)\frac{\partial f}{\partial \rho}+
D(0,t_1)\frac{\partial f}{\partial \alpha}dt_1+\\
\displaystyle
4\int\limits_{-\frac{\pi}{2}}^{\frac{\pi}{2}}2D(t_1,t_2)\frac{\partial f}{\partial \alpha}
dt_1dt_2 \biggr)\stackrel{f(\rho, \alpha)=f(\rho,2\pi-\alpha)}{=}
-4\int \limits_{S^1} D(0,t_1)\frac{\partial f}{\partial \alpha}dt_1+ 4\int\limits_{S^1}D(0,t_1)\frac{\partial f}{\partial \alpha}
dt_1=0.
\end{array}
$$
Therefore any point of the circle is a locally extremal point.
Hence the circle is locally extremal. The corollary is proved.

Now let us say a few worlds about M\"obius energy
which is (in the version from~\cite{Freed})
$$
f_M=\frac{1}{|\tau(t_1)-\tau(t_2)|^2}-\frac{1}{D^2(t_1,t_2)}.
$$
It has many remarkable properties (see~\cite{O-H1} and~\cite{Freed}).
M\"obius energies of homothetic knots are equal.
This energy is
invariant for M\"obius transformations (see also Section 5).
The variational equations and the gradient flow equation of M\"obius
energy was studied in~\cite{Freed}.

Unfortunately, for M\"obius energy, the variation $Var$ is always infinite,
and this mean that we can not perturb the knot in the way considered above.

The main property of M\"obius energy is as follows.
When a knot crossing tends to a double point, the energy tends to infinity.
The energy is always positive.
So every topological type of knot has a representative with minimal
value of energy, some normal form.

Notice that the main part of M\"obius energy is
$1/|\tau(t_1)-\tau(t_2)|^2$. The other part $1/D^2(t_1,t_2)$
is only a normalization that makes the integral convergent.
So let us make another normalization of the ``main part" of M\"obius
energy. In this case we often lose the invariance for M\"obius transformations.
Let us consider the following energy:
$$
\tilde f=\frac{D^3(x,y)}{|\tau(x),\tau(y)|^2}.
$$
It is easily seen that this energy on one hand has the above property and
on the other we can use our variational principles.
Note also that such an energy is the same for homothetic knots.
\begin{corollary}\label{c-3}
We present $V_1$ and $V_2$ for this energy:
$$
\begin{array}{l}
\begin{aligned}
\displaystyle
V_1(t_0)&=&\frac{2}{3R(t_0)}\Biggl(
4 \int \limits_{S^1} \Bigl(\frac{{|\tau(t)- \tau(t_0)|} ^3}{{D(t, t_0)} ^2}
\Bigl( 1-2\frac{R(t_0)}{{D(t,t_0)}}\Phi(t_0,t) \Bigr) \Bigr) dt-\\
\displaystyle
&&\frac{3}{\pi}\iint \limits_{S^1 \times S^1}\frac{{|\tau(t_2)- \tau(t_1)|} ^3}
{{D(t_2,t_1)} ^2}dt_1dt_2+
6\iint \limits_{A}
\frac{{|\tau(t_2)- \tau(t_1)|} ^2}{{D(t_2,t_1)} ^2} dt_1dt_2 \Biggr);
\end{aligned}\\
\displaystyle
V_2(t_0)=-\frac{8}{3R(t_0)} \int \limits_{S} \frac{{|\tau(t_1)- \tau(t_2)|} ^3}
{{D(t_0,t)} ^3}\Psi(t_0,t) dt.
\end{array}
$$
\end{corollary}\vspace{5mm}

\vspace{5mm}
\section{Definition and some basic properties of Mm-energy}

In this section we define the Mm-energy of a knot. The nature
of this energy differs from the energies considered in
the previous sections.

Let us fix some point $t_0$ on the circle and define
the real number $f_{Mm}(t_0)$.
Consider the map $\rho_{t_0}:S^1 \longrightarrow \r$ such that
$\rho_{t_0}(t)=|\tau (t)-\tau (t_0)|$. Let us note that
the map $\tau$ is smooth. Hence $\rho_{t_0}$ is also smooth
except for one point $t_0$.
If the number of maximums and minimums is finite, then we
define the function $f_{Mm}$ as follows:
$$
f_{Mm}(t_0)=\frac{1}{\rho_{t_0}(t_M)}+
\sum\limits_{t_{m_i}\in U_1} \frac{1}{\rho_{t_0}(t_{m_i})}-
\sum\limits_{t_{M_j}\in U_2} \frac{1}{\rho_{t_0}(t_{M_j})},
$$
where $t_M$ is one of the points where the function $\rho_{t_0}$
achieves its global maximum; $U_1$ is the set of all points of
the circle, except the point $t_0$, where the function
$\rho_{t_0}$ has local minimums; $U_2$ is the set of all points
of the circle, except the point $t_M$, where the function
$\rho_{t_0}$ has local maximums (see Fig. 5). Here we suppose
$t_0<t_*<t_0+2\pi$. In the case of an infinite number of maximums
and minimums we make a small smooth perturbation $\tilde
\rho_{t_0}$ so that the number of minimums and maximums becomes
finite. Now we can calculate the value of $\tilde f_{Mm}(t_0)$
for the function $\tilde \rho_{t_0}$ as it was made before.
Finally we define the $f_{Mm}(t_0)$ as the limit of $\tilde
f_{Mm}(t_0)$ in the $C^{\infty}$-topology.

\begin{figure}
$$\epsfbox{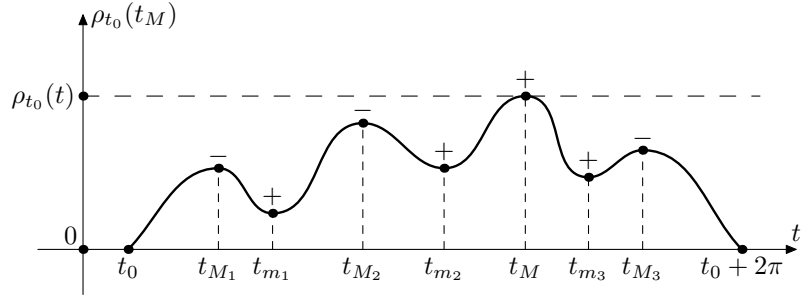}$$
\caption{The function $\rho_{t_0}$.}
\end{figure}

Now we define the Mm-energy.
\begin{definition}\label{d-3_1}
We call {\it Mm-energy} of the given knot the following number:
$$
E_{Mm}(\tau)=\int\limits_{S^1}f_{Mm}(t) dt,
$$
if the integral converges.
\end{definition}

\begin{remark}
Consider some small smooth perturbation of a knot.
Then for any point $t_0$ of the circle the function
$\rho_{t_0}$ is also perturbed in a smooth way.
At a generic point four possible modifications in
the sums of $f_{Mm}$ can occur: small changes of the values of
the maximums and minimums;
the death of one maximum and of the neighboring
minimum; conversely, the birth of one maximum and minimum at some
point; a local maximum close to the global maximum
can become the global maximum. In all these cases the variation
of the resulting $f_{Mm}$ is small.
This is the reason why the Mm-energy
depends on small perturbations of knots continuously.
\end{remark}

Further we formulate the basic properties of Mm-energy.

\begin{proposition}\label{p-3_1}
The Mm-energy is greater than or equal to $2$.
\end{proposition}

Consider the sum
$$
f_{Mm}(t_0)=\frac{1}{\rho_{t_0}(t_M)}+
\sum\limits_{t_{m_i}\in U_1} \frac{1}{\rho_{t_0}(t_{m_i})}-
\sum\limits_{t_{M_j}\in U_2} \frac{1}{\rho_{t_0}(t_{M_j})}
$$
We can fix the ordering of the minimums and the maximums in the standard
way:
$$
t_0<t_{M_1}<t_{m_1}<\ldots
<t_{M_k}<t_{m_k}<t_{M}<t_{m_{k+1}}<t_{M_{k+1}}<\ldots
<t_{m_{n}}<t_{M_{n}}<t_0+2\pi.
$$
Then we have
$$
\begin{array}{c}
\displaystyle
f_{Mm}(t_0)=
\sum \limits_{i=0}^{k}\Bigr(
\frac{1}{\rho_{t_0}(t_{m_i})}-\frac{1}{\rho_{t_0}(t_{M_i})} \Bigl)+
\frac{1}{\rho_{t_0}(t_M)}+
\sum \limits_{i=k+1}^{n}\Bigr(
\frac{1}{\rho_{t_0}(t_{m_i})}-\frac{1}{\rho_{t_0}(t_{M_i})} \Bigl)\ge\\
\displaystyle
0+\frac{1}{\rho_{t_0}(t_M)}+0=
\frac{1}{\rho_{t_0}(t_M)}.
\end{array}
$$
Finally, note that the length of the knot is $2\pi$, hence
the function $\rho_{t_0}(t_M)$ is smaller than or equal to $\pi$.
Therefore
$$
E_{Mm}(\tau)=\int\limits_{S^1}f_{Mm}(t) dt\le
\int\limits_{S^1}\frac{1}{\rho_{t}(t_{M_t})} dt\le
\int\limits_{S^1}\frac{1}{\pi} dt=\frac{2\pi}{\pi}=2.
$$
This completes the proof of Proposition~\ref{p-3_1}.

\begin{proposition}\label{p-3_2}
The Mm-energy is an invariant of homothety.
\end{proposition}

Suppose $\tau$ is a knot of length $2\pi$ and $\tilde \tau$ is
a homothetic knot of length $2l\pi$, where $l$ is the coefficient of
homothety. Then
$d\tilde t=ldt$ and $\tilde\rho(\tilde t)=l\rho(t)$ for any $t$,
and so $\tilde f_{Mm}(\tilde t)=f_{Mm}(t)/l$.
Thus we obtain
$$
E_{Mm}(\tilde\tau)=\int\limits_{S^1}\tilde f_{Mm}(\tilde t)d\tilde t=
\int\limits_{S^1}\frac{f_{Mm}(t)}{l}ldt=
\int\limits_{S^1}f_{Mm}(t)dt=E_{Mm}(\tau).
$$
Proposition~\ref{p-3_2} is proven.

So we can consider knots without any restriction on their lengths.

\begin{proposition}\label{p-3_3}
When two branches of the knot tends to a double crossing,
the Mm-energy tends to infinity.
\end{proposition}

Consider a smooth family $\{\tau_\lambda|\lambda\in[0,1]\}$ such that
$\tau_0$ is a smooth knot
with double crossing and $\tau_\lambda, \lambda\ne 0$ is a smooth knot without
any double crossing.
For every $\e$ we can choose a sufficiently small $\lambda$ satisfying the
following conditions:
there exist two points $t_1$ and $t_2$ with $|t_1-t_2|<\e^2$ such that
the functions
$\rho_{t_1}$ and $\rho_{t_2}$ have global minima at the points
$t_2$ and $t_1$ correspondingly; and the ball $B_{\e, p}$
of radius $\e$ with center
at the midpoint $p$ of the segment $[\tau_\lambda(t_1), \tau_{\lambda}(t_2)]$
has only two connected components of a knot $\tau_{\lambda}$ inside.

The family is smooth, hence the curvature of all knots is bounded by some $N$.
If $\e<1/N$, then
every point $t$ of the knot $\tau_{\lambda}$ inside the ball
$B_{\e /2,p}$ has one extremum (i.e., the
global minimum) of the function $\rho_t$ inside the ball $B_{\e ,p}$, and
every point $t$ of this knot inside the ball
$B_{\e, p}$ has no more than one extremum (i.e., the
global minimum) of $\rho _t$ inside the ball $B_{\e ,p}$.
Let us estimate the energy inside the ball $B_{\e,p}$.
$$
E_{Mm}(\tau_\lambda \cap B_{\e,p})>4\int\limits_{\frac{\e ^2}{2}}^{\frac{\e}{2}}
\frac{1}{t+\frac{\e ^2}{2}} dt=\left. 4\ln (t+ \frac{\e^2}{2})
\right|_{\frac{\e ^2}{2}}^{\frac{\e}{2}}=
4\ln \frac{\frac{e}{2}+\frac{e^2}{2}}{\e ^2}>4\ln \frac{2}{\e}.
$$
The other terms (we ignore the global minimum of $\rho_t$) of
the function $f_{Mm}$
changes in a smooth way, hence the Mm-energy grows to infinity.

Therefore Mm-energy separates knots from different topological
classes.

The following property is an essential property of Mm-energy.
\begin{proposition}\label{p-3_4}
The Mm-energy is well defined for piecewise smooth knots with obtuse angles.
\end{proposition}

If some point $t$ is ``near" the angle then the function $\rho_t$
is monotone function in some neighborhood of the vertex of an angle and
hence there are no minima or maxima of $\rho_t$ in this neighborhood.

In particular, the Mm-energy is well defined for piecewise linear knots
with obtuse angles.
So we can consider piecewise linear approximations of smooth
knots and take the restriction to the set of
piecewise linear knots. This property allows us to develop computer
experiments in calculating normal forms for Mm-energies
of topological classes of knots and the values of Mm-energies for this
normal forms.

Now we calculate Mm-energy for some knots.
First we find the Mm-energy of the circle $\tau_0$
$$
E_{Mm}(\tau_0)=\int\limits_{S^1}\frac{1}{2} dt=\pi.
$$
Unfortunately the circle is not the normal form for the class of
trivial knots. An example of the trivial knot with Mm-energy less
than $\pi$ is shown on Figure~6. This knot is a union of two arcs
of the circle. Direct calculations shows that the Mm-energy of
this knot is $2\ln(\frac{7+4\sqrt{3}}{3})\approx 3.070607<\pi$.

\begin{figure}
$$\epsfbox{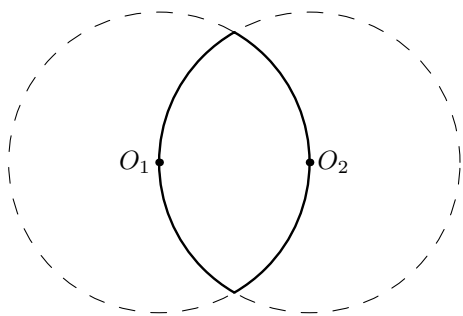}$$
\caption{Mm-energy of this knot is $2\ln(\frac{7+4\sqrt{3}}{3})$ .}
\end{figure}

Computer experiments provide upper bounds for the Mm-energies
of the normal forms for some topological classes (see the table behind).

\begin{center}
\vspace{0.2cm}
\begin{tabular}{|p{5cm}|c|}
\hline
{\sc CLASSES OF KNOTS} & {\sc the upper bounds for}\\
 & {\sc the energies of normal forms}\\
\hline
the class of the circle &  3.044012 \\
\hline
the class of the trefoil &  13.152759 \\
\hline
the class of the figure-eight &  19.450447 \\
\hline
the class of $5_1$ &  26.498108 \\
\hline
the class of $5_2$ &  27.168222 \\
\hline
the class of $6_1$ &  34.469191 \\
\hline
the class of $6_2$ &  35.466138 \\
\hline
the class of $6_3$ &  37.683129 \\
\hline
the class of the connected sum of right and left trefoils & 25.734616  \\
\hline
the class of the connected sum of two right trefoils & 26.748901  \\
\hline
\end{tabular}
\end{center}

\pagebreak

\vspace{0.2in}
\end{document}